between
can be
variable

%
\newlinechar`\^^J%
\message{^^J Attention, vous compilez un fichier comportant de nombreuses ^^J 
macros specifiques a l'Institut Fourier, il est recommande de lire ^^J
les commentaires figurant au debut de ce fichier. Ce sera  ^^J
necessaire si la compilation vous indique des erreurs.^^J^^J}
\message{ Please learn french and read the included comments in case of trouble.^^J
This file should compile with *plain* TeX. It inputs two files from the ^^J
AMS TeX distribution and requires AMS fonts: euler, symbol and extracm.^^J^^J}
%
%
\input amssym.def
\input amssym  
\catcode`\|=13 
\font\helvb=cmssbx10
\font\eightrm=cmr8
\font\eighti=cmmi8
\font\eightsy=cmsy8
\font\eightbf=cmbx8
\font\eighttt=cmtt8
\font\eightit=cmti8
\font\eightsl=cmsl8
\font\sixrm=cmr6
\font\sixi=cmmi6
\font\sixsy=cmsy6
\font\sixbf=cmbx6
\skewchar\eighti='177 \skewchar\sixi='177
\skewchar\eightsy='60 \skewchar\sixsy='60

\def\tenpoint{%
  \textfont0=\tenrm \scriptfont0=\sevenrm \scriptscriptfont0=\fiverm
  \def\rm{\fam0\tenrm}%
  \textfont1=\teni \scriptfont1=\seveni \scriptscriptfont1=\fivei
  \def\mit{\fam\@ne}\def\oldstyle{\fam1\teni}%
  \textfont2=\tensy \scriptfont2=\sevensy \scriptscriptfont2=\fivesy
    \def\itfam{4}\textfont\itfam=\tenit
  \def\it{\fam\itfam\tenit}%
  \def\slfam{5}\textfont\slfam=\tensl
  \def\sl{\fam\slfam\tensl}%
  \def\bffam{6}\textfont\bffam=\tenbf \scriptfont\bffam=\sevenbf
  \scriptscriptfont\bffam=\fivebf
  \def\bf{\fam\bffam\tenbf}%
  \def\ttfam{7}\textfont\ttfam=\tentt
  \def\tt{\fam\ttfam\tentt}%
  \abovedisplayskip=6pt plus 2pt minus 6pt
  \abovedisplayshortskip=0pt plus 3pt
  \belowdisplayskip=6pt plus 2pt minus 6pt
  \belowdisplayshortskip=7pt plus 3pt minus 4pt
  \smallskipamount=3pt plus 1pt minus 1pt
  \medskipamount=6pt plus 2pt minus 2pt
  \bigskipamount=12pt plus 4pt minus 4pt
  \normalbaselineskip=12pt
  \setbox\strutbox=\hbox{\vrule height8.5pt depth3.5pt width0pt}%
  \normalbaselines\rm}

\def\eightpoint{%
  \textfont0=\eightrm \scriptfont0=\sixrm \scriptscriptfont0=\fiverm
  \def\rm{\fam0\eightrm}%
  \textfont1=\eighti \scriptfont1=\sixi \scriptscriptfont1=\fivei
  \def\oldstyle{\fam1\eighti}%
  \textfont2=\eightsy \scriptfont2=\sixsy \scriptscriptfont2=\fivesy
  \textfont\slfam=\eightit
  \def\sl{\fam\itfam\eightit}%
  \textfont\slfam=\eightsl
  \def\sl{\fam\slfam\eightsl}%
  \textfont\bffam=\eightbf \scriptfont\bffam=\sixbf
  \scriptscriptfont\bffam=\fivebf
  \def\bf{\fam\bffam\eightbf}%
  \textfont\ttfam=\eighttt
  \def\tt{\fam\ttfam\eighttt}%
  \abovedisplayskip=9pt plus 2pt minus 6pt
  \abovedisplayshortskip=0pt plus 2pt
  \belowdisplayskip=9pt plus 2pt minus 6pt
  \belowdisplayshortskip=5pt plus 2pt minus 3pt
  \smallskipamount=2pt plus 1pt minus 1pt
  \medskipamount=4pt plus 2pt minus 1pt
  \bigskipamount=9pt plus 3pt minus 3pt
  \normalbaselineskip=9pt
  \setbox\strutbox=\hbox{\vrule height7pt depth2pt width0pt}%
  \normalbaselines\rm}

\font\petcap=cmcsc10

\tenpoint
\hsize=12.5cm
\vsize=19cm
\parskip 5pt plus 1pt
\parindent=1cm
\baselineskip=13pt
\hoffset=-0.1cm 
\def\footnoterule{\kern-6pt
  \hrule width 2truein \kern 5.6pt} 

\def\ie{{\sl i.e.\ }}

\def\omini{\raise 1ex\hbox{\ept o}}
\def\emini{\raise 1ex\hbox{\ept e}}
\def\ermini{\raise 1ex\hbox{\ept er}}
\def\remini{\raise 1ex\hbox{\ept re}}

\def\lead{\leaders\hbox to 10pt{\hss.\hss}\hfill}
\def\somt#1|#2|{\vskip 8pt plus1pt minus 1pt
                \line{#1\lead #2}}
\def\soms#1|#2|{\vskip 2pt
                \line{\qquad #1\lead #2}}  
\def\somss#1|#2|{\vskip 1pt   
                 \line{\qquad\qquad #1\lead #2}} 


\def\aujour{\ifnum\day=1 1\ermini\else\number\day\fi\
\ifcase\month\or janvier\or f\'evrier\or mars\or avril\or mai\or juin\or
juillet\or aout\or septembre\or octobre\or novembre\or d\'ecembre\fi\
\number\year}
\def\today{\ifcase\month\or january \or february \or march \or april
\or may \or june\or july\or august \or september\or october\or november\or
december\fi\ \number\day , \number\year}

\newskip\afterskip
\catcode`\@=11
\def\p@int{.\par\vskip\afterskip\penalty100} 
\def\p@intir{\discretionary{.}{}{.\kern.35em---\kern.7em}}
\def\pointir{\afterassignment\pointir@\global\let\next=}
\def\pointir@{\ifx\next\par\p@int\else\p@intir\fi\egroup\next}
\catcode`\@=12
\def|{\relax\ifmmode\vert\else\findef\fi}
\def\findef{\errhelp{Cette barre verticale ne correspond ni a un \vert
mathematique
                        ni a une fin de definition, le contexte doit vous
indiquer ce qui manque.
                        Si vous vouliez inserer un long tiret, le codage
recommande est ---,
                        dans tous les cas, la barre fautive a ete supprimee.}%
                        \errmessage{Une barre verticale a ete trouvee en
mode texte}}

\def\TITR#1|{\null{\mss\baselineskip=17pt
                           \vskip 3.25ex plus 1ex minus .2ex
                           \leftskip=0pt plus \hsize
                           \rightskip=\leftskip
                           \parfillskip=0pt
                           \noindent #1
                           \par\vskip 2.3ex plus .2ex}}
 
\def\auteur#1|{\penalty 500
               \vbox{\centerline{\si
                 \iffrance par \else by \fi #1}
                \vskip 10pt}\penalty 500}

\def\resume#1|{\penalty 100
                           {\leftskip=\parindent
                            \rightskip=\leftskip
                            \eightpoint\bgroup\petcap \skip\afterskip=0pt
                             \iffrance R\'esum\'e \else Abstract \fi\pointir
                            #1 \par}
                           \penalty -100}

\def\titre#1|{\null\baselineskip14pt
                           {\helvb
                           \vskip 3.25ex plus 1ex minus .2ex
                           \leftskip=0pt plus \hsize
                           \rightskip=\leftskip
                           \parfillskip=0pt
                           \noindent #1
                           \par\vskip 2.3ex plus .2ex}}


\def\section#1|{
                                \bgroup\bf
                                 \par\penalty -500
                                 \vskip 3.25ex plus 1ex minus .2ex
                                 \skip\afterskip=1.5ex plus .2ex
                                  #1\pointir}

\def\ssection#1|{
                                 \bgroup\petcap
                                  \par\penalty -200
                                  \vskip 3.25ex plus 1ex minus .2ex
                          \skip\afterskip=1.5ex plus .2ex
                                   #1\pointir}
 

\def\th#1|{ 
                   \bgroup \sl
                        \def\findef{\egroup\par}
                        \bgroup\petcap
                         \par\vskip 2ex plus 1ex minus .2ex
                         \skip\afterskip=0pt
                           #1\pointir}

\def\defi#1|{ 
                   \bgroup \rm
                        \def\findef{\egroup\par}
                        \bgroup\petcap
                         \par\vskip 2ex plus 1ex minus .2ex
                         \skip\afterskip=0pt
                           #1\pointir}

\def\rque#1|{\bgroup \sl
                          \par\vskip 2ex plus 1ex minus .2ex\skip\afterskip=0pt
                          #1\pointir}

\def\dem{\bgroup \sl
                  \par\vskip 2ex plus 1ex minus .2ex\skip\afterskip=0pt
                  \iffrance D\'emonstration\else Proof\fi\pointir}

\def\preuve{\bgroup \sl
                  \par\vskip 2ex plus 1ex minus .2ex\skip\afterskip=0pt
                  \iffrance Preuve\else Proof\fi\pointir}


\def\_#1{_{\baselineskip=.7 \baselineskip
                                                                       
\vtop{\halign{\hfil$\scriptstyle{##}$\hfil\cr #1\crcr}}}}

\def\build#1#2\fin{\mathrel{\mathop{\kern0pt#1}\limits#2}}

\def\frac#1/#2{\leavevmode\kern.1em
   \raise.5ex\hbox{$\scriptstyle #1$}\kern-.1em
      /\kern-.15em\lower.25ex\hbox{$\scriptstyle #2$}}

{\obeylines
\gdef\iffin{\parskip0pt\parindent0pt 
            \vskip1cm
            \centerline{--$~\diamondsuit~$--}
            \vskip1cm
            \eightpoint
            Universit\'e de Grenoble I
            {\bf Institut Fourier}
            Laboratoire de Math\'ematiques
            associ\'e au CNRS (URA 188)
            B.P. 74
            38402 ST MARTIN D'H\`ERES Cedex (France)
            \vskip1cm
            (\aujour)
            }}

\newif\ifchiffre
\def\chiffre{\chiffretrue}
\chiffre
\newdimen\laenge
\def\lettre#1|{\setbox3=\hbox{#1}\laenge=\wd3\advance\laenge by 3mm
\chiffrefalse}
\def\article#1|#2|#3|#4|#5|#6|#7|%
    {{\ifchiffre\leftskip=7mm\noindent
     \hangindent=2mm\hangafter=1
\llap{[#1]\hskip1.35em}\bgroup\petcap #2\pointir {\sl #3}, {\rm #4}
\nobreak{\bf #5}
({\oldstyle #6}), \nobreak #7.\par\else\noindent \advance\laenge by 4mm
\hangindent=\laenge\advance\laenge by -4mm\hangafter=1
\rlap{[#1]}\hskip\laenge\bgroup\petcap #2\pointir {\sl #3}, #4 {\bf #5}
({\oldstyle
#6}), #7.\par\fi}} 
\def\livre#1|#2|#3|#4|#5|%
    {{\ifchiffre\leftskip=7mm\noindent
    \hangindent=2mm\hangafter=1
\llap{[#1]\hskip1.35em}\bgroup\petcap #2\pointir{\sl #3}, #4, {\oldstyle
#5}.\par
\else\noindent
\advance\laenge by 4mm \hangindent=\laenge\advance\laenge by -4mm
\hangafter=1
\rlap{[#1]}\hskip\laenge\bgroup\petcap #2\pointir
{\sl  #3}, #4, {\oldstyle #5}.\par\fi}}
\def\divers#1|#2|#3|#4|%
    {{\ifchiffre\leftskip=7mm\noindent
    \hangindent=2mm\hangafter=1
     \llap{[#1]\hskip1.35em}\bgroup\petcap #2\pointir #3, {\oldstyle #4}.\par
\else\noindent
\advance\laenge by 4mm \hangindent=\laenge\advance\laenge by -4mm
\hangafter=1
\rlap{[#1]}\hskip\laenge\bgroup\petcap #2\pointir #3,{\oldstyle #4}.\par\fi}}
\def\div#1|#2|#3|#4|
{{\ifchiffre\leftskip=7mm\noindent
\hangindent=2mm\hangafter=1
\llap{[#1]\hskip1.35em}\bgroup\petcap #2\pointir {\sl  #3},{\oldstyle #4}.\par
\else\noindent
\advance\laenge by 4mm \hangindent=\laenge\advance\laenge by -4mm
\hangafter=1
\rlap{[#1]}\hskip\laenge\bgroup\petcap #2\pointir {\sl  #3},{\oldstyle
#4}.\par\fi}}


\font\si=cmssi10

\font\mss=cmss12 scaled \magstep1 

\hsize=12.5cm
\vsize=19cm
\parskip 5pt plus 1pt
\parindent=1cm
\baselineskip=13pt

  \hoffset=0.7cm
 \voffset=0.8cm

\font\cms=cmbsy10 at 5pt
\font\gros=cmex10 at 17.28pt

\newif\iffrance

\def\anglais{\francefalse}

\def\TITR#1|{\null{\mss\baselineskip=17pt
                           \vskip 3.25ex plus 1ex minus .2ex
                           \leftskip=0pt plus \hsize
                           \rightskip=\leftskip
                           \parfillskip=0pt
                           \noindent #1
                           \par\vskip 2.3ex plus .2ex}}
 
\def\auteur#1|{\penalty 500
               \vbox{\centerline{\si
                 \iffrance par \else by \fi #1}
                \vskip 10pt}\penalty 500}

\def\resume#1|{\penalty 100
                           {\leftskip=\parindent
                            \rightskip=\leftskip
                            \eightpoint\bgroup\petcap \skip\afterskip=0pt
                             \iffrance R\'esum\'e \else Abstract \fi\pointir
                            #1 \par}
                           \penalty -100}

\def\titre#1|{\null\penalty-500\baselineskip14pt
                           {\helvb
                           \vskip 3.25ex plus 1ex minus .2ex
                           \leftskip=0pt plus \hsize
                           \rightskip=\leftskip
                           \parfillskip=0pt
                           \noindent #1
                           \par\nobreak\vskip 2.3ex plus .2ex}\penalty 5000}


\def\summ{ \mathop{\sum}\limits}
\def\cupp{\mathop{\bigcup}\limits}

\def\hfl{{\hbox to 12mm{\rightarrowfill}}}

\def\egal{{{\lower5pt\hbox to 0.5cm{$\hrulefill$}}\atop 
{\raise10pt\hbox to 0.5cm{$\hrulefill$}}}}

\def\upuparrow{\eqalign{&\big\uparrow\cr\noalign{\vskip-14pt}&\kern-1.8pt
               \uparrow\cr}}
\def\ddarrow{\eqalign{&\raise15pt\hbox{$\big\downarrow$}\cr
             \noalign{\vskip-25pt} &\big\downarrow\cr}}

\def\cu{\hbox{\cms\char'133}}
\def\ca{\hbox{\cms\char'134}}  
 \def\hookup{{\lower6.3pt\hbox{\cu}\kern-4.32pt\big\uparrow}}
 \def\hookdown{\raise5.3pt\hbox{\ca}\kern-4.32pt\lower3pt\hbox{\big\downarrow}}

\def\grtilde{\kern10pt\hbox{\gros\char'147}}

\let\le=\leq 

\def\S{\mathhexbox278\kern.15em}
\def\\ {\smallsetminus}

\def\ept{\eightpoint}

\def\ld {, \ldots,}

\def\coin{\mathrel{\raise2pt\hbox{$\scriptstyle
|$}\kern-2pt\lower2pt\hbox{$-$}}}
\def\rect{\mathrel{\lower2pt\hbox{$-$}\kern-2pt\raise2pt\hbox{$\scriptstyle
|$}}} 


\def\im{\mathop{\rm Im}\nolimits}

{\obeylines
\gdef\iffin{\parskip0pt\parindent0pt 
            \vskip1cm
            \centerline{--$~\diamondsuit~$--}
            \vskip1cm
            \eightpoint
            Universit\'e de Grenoble I
            {\bf Institut Fourier}
            Laboratoire de Math\'ematiques
            associ\'e au CNRS (URA 188)
            B.P. 74
            38402 ST MARTIN D'H\`ERES Cedex (France)
            \vskip1cm
            (\aujour)
            }}


\long\def\nomm#1|#2|#3|{\line{$\vtop{\hsize2.5cm #1~}\vtop{\hsize1cm #2~}
\vtop{\hsize=9cm\normalbaselines\parshape 1 0cm 9cm #3.}$}
\medskip}

\long\def\nom#1|#2|{\centerline{$\vtop{\hsize=3cm{\bf #1}~: }
\vtop{\hsize=9.5cm\normalbaselines\parshape 1 0cm 9.5cm #2.}$}
\medskip}

\def\boxit#1#2{\hbox{\vrule
 \vbox{\hrule\kern#1
  \vtop{\hbox{\kern#1 #2\kern#1}%
   \kern#1\hrule}}%
 \vrule}}

\newbox\texte
\def\texteencadre#1|#2|#3|{\setbox1=\vbox{#3}
\setbox\texte=\vbox{\hrule height#1pt%
\hbox{\vrule width#1pt\kern#2pt\vbox{\kern#2pt \hbox{\box1}\kern#2pt}%
\kern#2pt\vrule width#1pt}\hrule height#1pt}}

\def\blanc{\hbox to 0pt{\vrule height13pt depth5pt width0pt}}
\def\blancs{\hbox to 0pt{\vrule height9pt depth5pt width0pt}}
\def\blan{\hbox to 0pt{\vrule height6pt depth5pt width0pt}}
\def\bla{\hbox to 0pt{\vrule height7pt depth5pt width0pt}}
\def\hrulefill{\leaders\hrule height0.2pt\hfill}

\def\bb#1&#2\cr{
\hbox to 1cm{\hfil\strut#1~\vrule}
\hbox to 1cm{\hfil\strut#2~\vrule}
}
\def\aa#1&#2&#3&#4&#5&#6&#7&#8&#9\cr{\line{$
\hbox to 3cm{\vrule width 1pt \strut~#1\hfil\vrule}
\hbox to 1cm{\hfil\strut#2~\vrule}
\hbox to 1cm{\hfil\strut#3~\vrule}
\hbox to 1cm{\hfil\strut#4~\vrule}
\hbox to 1cm{\hfil\strut#5~\vrule}
\hbox to 1cm{\hfil\strut#6~\vrule}
\hbox to 1cm{\hfil\strut#7~\vrule}
\hbox to 1cm{\hfil\strut#8~\vrule width 1pt}
\bb#9\cr
$\hfil}}


           \def\rb{{\Bbb R}}


\font\tenmib=cmmib10 
\font\sevenmib=cmmib7
\font\fivemib=cmmib5
\expandafter\chardef\csname pre boldmath.tex at\endcsname=\the\catcode`\@
\catcode`\@=11
\def\hexanumber@#1{\ifcase#1 0\or 1\or 2\or 3\or 4\or 5\or 6\or 7\or 8\or
 9\or A\or B\or C\or D\or E\or F\fi}

\skewchar\tenmib='177 \skewchar\sevenmib='177 \skewchar\fivemib='177

\newfam\mibfam   
\textfont\mibfam=\tenmib \scriptfont\mibfam=\sevenmib
\scriptscriptfont\mibfam=\fivemib
\def\mib@hex{\hexanumber@\mibfam}
\mathchardef\bfalpha="0\mib@hex 0B
\mathchardef\bfbeta="0\mib@hex 0C
\mathchardef\bfgamma="0\mib@hex 0D
\mathchardef\bfdelta="0\mib@hex 0E
\mathchardef\bfepsilon="0\mib@hex 0F
\mathchardef\bfzeta="0\mib@hex 10
\mathchardef\bfeta="0\mib@hex 11
\mathchardef\bftheta="0\mib@hex 12
\mathchardef\bfiota="0\mib@hex 13
\mathchardef\bfkappa="0\mib@hex 14
\mathchardef\bflambda="0\mib@hex 15
\mathchardef\bfmu="0\mib@hex 16
\mathchardef\bfnu="0\mib@hex 17
\mathchardef\bfxi="0\mib@hex 18
\mathchardef\bfpi="0\mib@hex 19
\mathchardef\bfrho="0\mib@hex 1A
\mathchardef\bfsigma="0\mib@hex 1B
\mathchardef\bftau="0\mib@hex 1C
\mathchardef\bfupsilon="0\mib@hex 1D
\mathchardef\bfphi="0\mib@hex 1E
\mathchardef\bfchi="0\mib@hex 1F
\mathchardef\bfpsi="0\mib@hex 20
\mathchardef\bfomega="0\mib@hex 21
\mathchardef\bfvarepsilon="0\mib@hex 22
\mathchardef\bfvartheta="0\mib@hex 23
\mathchardef\bfvarpi="0\mib@hex 24
\mathchardef\bfvarrho="0\mib@hex 25
\mathchardef\bfvarsigma="0\mib@hex 26
\mathchardef\bfvarphi="0\mib@hex 27
\mathchardef\bfimath="0\mib@hex 7B
\mathchardef\bfjmath="0\mib@hex 7C
\mathchardef\bfell="0\mib@hex 60
\mathchardef\bfwp="0\mib@hex 7D
\mathchardef\bfpartial="0\mib@hex 40
\mathchardef\bfflat="0\mib@hex 5B
\mathchardef\bfnatural="0\mib@hex 5C
\mathchardef\bfsharp="0\mib@hex 5D
\mathchardef\bftriangleleft="2\mib@hex 2F
\mathchardef\bftriangleright="2\mib@hex 2E
\mathchardef\bfstar="2\mib@hex 3F
\mathchardef\bfsmile="3\mib@hex 5E
\mathchardef\bffrown="3\mib@hex 5F
\mathchardef\bfleftharpoonup="3\mib@hex 28
\mathchardef\bfleftharpoondown="3\mib@hex 29
\mathchardef\bfrightharpoonup="3\mib@hex 2A
\mathchardef\bfrightharpoondown="3\mib@hex 2B
\mathchardef\bflhook="3\mib@hex 2C 
\mathchardef\bfrhook="3\mib@hex 2D 
\mathchardef\bfldotp="6\mib@hex 3A 
\catcode`\@=\csname pre boldmath.tex at\endcsname
\anglais

\def\oo{\overline\omega }
\def\Herm{\mathop{\rm Herm}\nolimits}
\def\im{\mathop{\rm Im}\nolimits}

\TITR A REPRESENTATION OF ISOMETRIES ON FUNCTION SPACES|

\auteur Mikhail G. ZAIDENBERG|

\centerline {\bf INTRODUCTION}

This paper contains a proof of the main result previously announced in
[Za1] and of its generalization to a class of ideal spaces.  Namely,
we prove the validity of a weighted shift representation of the
surjective isometries between Banach function spaces which satisfy
some minor restrictions. This class encludes at least all
rearrangement-invariant (r.i.)  spaces. Our proof follows Lumer's
scheme [Lu] and uses some ideas due to A.  Pelczynski (see [Ro];
cf. also [BrSe] and [SkZa1,2]). Recall that in [Za1] all the spaces
are considered over the field $\bf C$ of complex numbers.

Due to a recent revival of interest in the isometric theory of
Banach function  spaces  I have been asked several times during the last
few years by my colleagues working on the subject about the proofs of the
results announced in [Za1]. As a matter of fact, the  proofs of Theorem 1
and Proposition 1 from [Za1]  were published in [Za2], a  Russian journal
with a rather limited  circulation.  An English translation of that
journal though prepared has never appeared  due to some circumstances.
Nevertheless, a translation of my article
[Za2] was circulating among a small number of
experts (see e.g. references in [KaRa2]).
In order to satisfy numerous requests of my colleagues  I am reproducing
here this translation, certainly updaiting and modifying it.

Let me briefly mention several further developments.  A generalization
of Theorem 1 in [Za1] to the real case has been recently obtained in
[KaRa1,2]. Papers [KaRa1,2] contain also a new proof of Theorem 4 in
[Za1], which yields a characterization of the $L_p$-spaces as
r.i. spaces with non-standard isometry groups (my original proof of
this theorem covered both the real and complex cases but it was never
published). Some corollaries of the main results of [Za1] are extended
to the real case by different methods in [AbZa].  New proofs of the
remaining statements from [Za1] can be found in the forthcoming papers
[PKL], [Ra].  Some additional information can be also found in [FlJe].

\bigskip

\centerline{\bf HERMITIAN OPERATORS} 

A Banach space $E$ of measurable functions on a measure space $(\Omega
,\Sigma ,\mu )$ is called {\it an ideal space} if $f \in E,~~|g|\le
|f|$ and $g\in L^0(\Omega ,\Sigma ,\mu )$ imply that $g\in E$ and
$\|g\| \le \|f\|$.

If additionally, the  equimeasurability of functions $|f|$ and $|g|$,
where $f\in E$ and $g\in L^0$, implies that
$g\in E$ and $\|g\| = \|f\|$, then the space $E$ is called {\it symmetric}
or rearrangement-invariant.

Let $\Sigma_0 = \{ \sigma \in \Sigma \mid \mu (\sigma )<\infty\}$.  We
will always assume that the characteristic function $\chi _\sigma \in
E$ of every set $\sigma \in \Sigma _0$ belongs to $E$. Let
$\rho_\sigma $ be the projector $\rho _\sigma (x) = \chi _\sigma x$,
$x\in E$. The image $\rho _\sigma E$ is called {\it a band}, or {\it a
component} of the space $E$.

\defi Definition 1|An ideal space $E$ is called {\sl projection  
provided}, if for each finite collection $\oo  = \{\omega
_i\}^n_{i=1}$ of disjoint sets $\omega _i \in \Sigma _0$, $i=1\ld n$,
there exists a projector $\rho ^c_{\oo }$ of norm 1
on the subspace of step-functions 
$$E^n(\oo) = \Bigl\{\summ^n_{i=1}
c_i\chi _{\omega _i}, ~c_i \in {\bf C},~ i=1\ld n\Bigr\}$$ that commutes
with the projectors $\rho _{\omega _i}$, $i=1\ld n$.|

Any symmetric space is projection provided; we can take for $\rho
^c_{\oo}$ the Haar's projector of the conditional expectation ([SMB],
p.~95).

\defi Definition 2|We will say that the space $E$ is {\it free from
$L_2$-components} if for each component the norm of $E$ does not
coincide with the norm  of $L_2(\nu )$, where $\nu $ is a positive
measure on $\Sigma $.

An operator $H$ on $E$ is called {\it Hermitian} if $\|e^{irH}\| =1$,
$\forall t\in \rb$. The set of Hermitian operators is denoted as $\Herm
(E)$. By $LM^r_\infty(E)$ we denote the subset of the operators of
multiplication by bounded real functions.|

For simplicity we assume in the paper that $\Omega$ is $[0,1]$ or a
line with the Lebesgue measure $\mu$.

\th Proposition 1|Let $E$ be a projection provided Banach ideal
space. If $E$ is free from $L_2$-components, then $\Herm(E) =
LM^r_\infty(E)$.|

\preuve It is enough to prove that any Hermitian operator $H$ on $E$
holds the following property: $$\chi _{\Omega \setminus \omega } \cdot
H \chi _\omega = 0,~~ \forall \omega \in \Sigma _0.\eqno(1)$$ Indeed,
if (1) is true, then $$\chi _{\omega _1}\cdot H\chi _{\omega _2} =
\chi _{\omega _2}\cdot H\chi _{\omega _1} ~ (\, =H \chi _{\omega
_1\cap \omega _2}),~~ \forall \omega_1, \omega_2 \in \Sigma
_0.\eqno(2)$$ So, the equality $$\chi _\omega h = H\chi _\omega
,~~\omega \in \Sigma _0\eqno(3)$$ determines a measurable function $h$
on $\Omega $ such that $$e^{itH}(\chi _\omega ) = \chi _\omega
e^{ith}.\eqno(4)$$ By Theorem 2 of [Za3], (4) implies that
$$e^{itH}(f) = e^{ith}\cdot f,~~\forall f\in E,\eqno(5)$$ and
therefore $$H(f) = h\cdot f,~~\forall f\in E,\eqno(6)$$ and $\im h=0$,
$\|H\| = \|h\|_{L_\infty}$ [Lu, Lemma 7].

Suppose that there exists an operator $H_0 \in \Herm(E)$ which does
not hold (1), that means that $\chi _{\Omega \setminus \omega _0} H_0
\chi _{\omega _0} \ne 0$ for some $\omega _0 \in \Sigma _0$. Choose
disjoint sets $\omega _1\ld \omega _n \in \Sigma _0$ ($\omega _i\cap
\omega _j = \emptyset$, $i\ne j,\,i,j=0,\dots,n$) such that $$|\chi
_{\omega _i} H_0 \chi _{\omega _0} - \lambda _i\chi _{\omega _i}| <
|\lambda _i|,\eqno(7)$$ where $\lambda _i\ne 0$, $i=1,\dots,n$. Let
$\oo = \{\omega _i\}^n_{i=0}$. Due to
 Lumer's Lemma [Lu, Lemma 8], the operator $H^{\oo}_0$ determined by
the equality $$H^{\oo}_0 = \rho ^c_{\oo} \cdot H_0 \rho ^c_{\oo}$$ is
a Hermitian operator on the subspace $E^{n+1}(\oo)$. Since $E$ is an
ideal space, $\|\rho ^c_{\oo}\| =1$, and $\rho ^c_{\oo}(\chi _{\omega
_i}x) = \chi _{\omega _i} \rho ^c_{\oo}(x)$, $i=0,1\ld n$, from (7) it
follows that $$\eqalign{ \|\chi _{\omega _i} H^{\oo}_0 \chi _{\omega
_0} - \lambda _i\chi _{\omega _i}\| &= \|\rho ^c_{\oo} (\chi _{\omega
_i} H_0\chi _{\omega _0} - \lambda _i\chi _{\omega _i})\|\cr &\le
\|\chi _{\omega _i} H_0\chi _{\omega _0}-\lambda _i\chi _{\omega
_i}\|\cr & < |\lambda _i|\cdot \|\chi _{\omega _i}\|,~~ i=1\ld
n,\cr}$$ and whence $$\chi _{\omega _i} H^{\oo}_0 \chi _{\omega _0}
\ne 0,~~ i=1\ld n.\eqno(8)$$

Next we show that the subspace $E^{n+1}_{\oo}$ is Euclidean, and
$$\Bigl\|\summ^n_{i=0} c_i\chi _{\omega _i}\Bigr\|^2 = \summ^n_{i=0}
|c_i|^2 \|\chi _{\omega _i}\|^2.\eqno(9)$$ To this point we use the
following lemma \footnote{(*)}{Similar statements can be found, for
instance, in [KaWo] (complex case), [SkZa1,2] (real case); see also
the bibliography therein.}.

\th Lemma 1|Let $E^{n+1}$ be an ideal Minkowski space
\footnote{(**)}{\rm i.e. a finite dimensional Banach space} over $\bf
C$, $\{e_i\}^n_0$ be the standard basis in $E^{n+1}$, and $H$ be a
Hermitian operator on $E^{n+1}$ such that $$(He_0,e_k) \ne 0,~~k=1\ld
n.\eqno(8')$$ Then $E^{n+1}$ is a Euclidean space and
$$\Bigl\|\summ^n_{i=0} c_ie_i\Bigr\|^2 = \summ^n_{i=0} |c_i|^2
\|e_i\|^2.\eqno(9')$$|

\preuve Let $G_0$ be the connected component of unity in the isometry
group ${\rm Iso}\,(E^{n+1})$, and let $\langle \cdot,\cdot\rangle $ be
a $G_0$-invariant scalar product such that $\langle e_0,e_0\rangle
=1$. The orbit $G_0e_0$ is a connected $G_0$-invariant submanifold in
$E^{n+1}$, which is contained in the intersection of the Minkowski
sphere $S(E^{n+1})$ and the Euclidean sphere $S^{2n+1}$ (indeed, $G_0
\subset U(n)$ by our choice of scalar product). The tangent space $T$
to the orbit $G_0e_0$ at point $e_0$ is invariant with respect to the
stationary subgroup $G_0^{e_0}$ of $e_0$. Note that the operator
$t_k(\varphi )$ of the rotation of $k$--th coordinate on angle
$\varphi $ (i.e. $t_k(\varphi ) e_k = e^{i\varphi }e_k,~ t_k(\varphi
)e_i = e_i,~i\ne k)$ belongs to the subgroup $G^{e_0}_0$ and,
therefore, $t_k(\varphi ) T=T$. The latter is possible either if the
$k$--th coordinate of any vector of $T$ is equal to zero, or if
$e_k\in T$ and $ie_k\in T$. Since the tangent vector $iHe_0$ to the
curve $e^{itH}(e_0)\subset G_0e_0$ belongs to the subspace $T$, in
view of $(8')$, the $k$--th coordinate of this vector is
non--zero. Therefore, $e_k\in T$ and $ie_k\in T$ $(k=1\ld n)$. Besides
that, the subspace $T$ contains vector $ie_0$ tangent to the curve
$e^{itH}(e_0) \subset G_0e_0$. Thus, the real dimension of the
subspace $T$ and, hence, also of the orbit $G_0e_0$ is equal to
$2n+1$. Since $G_0e_0$ is a compact connected manifold, we have
$$G_0e_0 = S(E^{n+1}) = S^{2n+1}~,\eqno(10)$$ that proves that
$E^{n+1}$ is Euclidean.

Because eigen-vectors $e_k$ and $e_\ell$ $(\ell \ne k)$ of the
operator $t_k(\pi ) \in G_0$ are orthogonal, $\{e_k\}^n_0$ is an
orthogonal basis in $E^{n+1}$. This implies (9'). The lemma is proved.

Let us come back to the proof of Proposition 1. Let $S = \cupp^n_{i=1}
\omega _i$, $E^b(S)$ be the closure in $\rho _S E$ of the set of
finite-valued (step) functions. Passing to the limit over
subpartitions, we can prove that $E^b(S)$ is a Hilbert space. Due to
reflexivity of the space $E^b(S)$, the norm on $E^b(S)$ is absolutely
continuous [Zab, Theorem 30], \ie $\|\chi _\sigma \| \to 0$ as $\mu
(\sigma ) \to 0$, and $E^b(S) = (E^b(S))''$. Therefore, $E^b(S) = \rho
_SE$. Absolute continuity of the norm permits to get the equality
$$\|\chi _\sigma \|^2 = \summ^\infty_{i=1} \|\chi _{\sigma
_i}\|^2,\eqno(11)$$ where $\sigma _i$, $\sigma \in \Sigma (S) = \Sigma
\cap S$, $\sigma _i \cap \sigma _j = \emptyset$, $i\ne j$, $\sigma =
\cupp^\infty_{i=1} \sigma _i$.

Let $\nu (\sigma ) = \|\chi _\sigma \|^2$. Due to (11), $\nu $ is a
positive (absolutely continuous with respect to $\mu $)
$\sigma$--additive measure on the algebra $\Sigma (S)$. Previous
arguments show that $\rho _SE = L_2(\nu )$, that contradicts to our
assumption. The proposition is proved.

\th Proposition 2|Let $E$ be a symmetric space such that the norm on
$E$ is not proportional to the norm of the space $L_2(\Omega ,\Sigma
,\mu )$. Then $\Herm(E) = LM^r_\infty(E)$.|

The proof is quite similar to that of Proposition 1. A partition $\oo$
is chosen in a way that provides (8) just for $i=1$, but we should
require, besides that, that $\mu (\omega _i) = \mu (\omega _1)$,
$i=2\ld n$. Using Lemma IX.8.4 of [Ro] instead of Lemma 1, we obtain
the equality $$\Bigl\| \summ^n_{i=1} c_i\chi _{\omega _i}\Bigr\|^2 =
\|\chi _{\omega _1}\|^2 \summ^n_{i=1} |c_i|^2.\eqno(9'')$$ Passing to
the limit over partitions we can see that $\|\cdot \|_E =
k\|\cdot\|_{L_2(\mu )}$ on any component $\rho _s E$, $s\in \Sigma
_0$, where $k = \varphi _E(1)$ (here $\varphi _E$ denotes the
fundamental function of the symmetric space $E$, i.e. $\varphi _E(t) =
\|\chi _\sigma \|$, where $\mu (\sigma )=t$). This yields a
contradiction and completes the proof.

\bigskip

\centerline{\bf MAIN THEOREM} 

\th Theorem 1|

(a) Let $E_i = E_i(\Omega _i)$ $(i=1,2)$ be projection provided
 ideal Banach spaces, free from $L_2$-components. Then for any
isometric isomorphism $Q: E_1 \to E_2$ there exists a measurable
function $q$ and an invertible measurable transformation $\varphi
:\Omega _2 \to \Omega_1 $ such that $$(Qf)(t) = q(t) f(\varphi
(t)),~~\forall f\in E_1.\eqno(12)$$

(b) The same conclusion is true providing that $E_i$ $(i=1,2)$ are
symmetric spaces and the norm on $E_1$ is not proportional to the norm
on $L_2(\mu )$.|

\preuve As it was proved in [Se], the norm of a Hilbert symmetric
space is proportional to $L_2$-norm. Therefore, under the conditions
in {\sl (b)}, $\|\cdot \|_{E_2} \ne k\| \cdot \|_{L_2(\mu )}$. By
Propositions 1 and 2, $\Herm(E_i) = LM^r_\infty(E_i)$,
$i=1,2$. Consider the mapping $Q_* : \Herm(E_1) \to \Herm (E_2)$,
$Q_*(H) = QHQ^{-1}$. It is not difficult to verify that $Q_*$ is an
algebraic isomorphism. It generates an isomorphism of Boolean algebras
$\theta _* : \Sigma _1 \to \Sigma _2$, connected with $Q$ by the
relation $Q_*(\chi _\sigma ) = \chi _{\theta _* \sigma }$, $\sigma \in
\Sigma _1$. Define a measurable transformation $\varphi : \Omega _2
\to \Omega _1$ by the equalities $$\chi _{\theta _*\sigma }\varphi =
Q_*(t\chi _\sigma ),~~\sigma \in (\Sigma _1)_0.\eqno(13)$$ Obviously,
$$Q_*f = f(\varphi ) ~ \forall f\in L^r_\infty(\Omega _1),\quad \chi
_{\theta _*\sigma } = \chi _\sigma (\varphi ),~\sigma \in \Sigma
_1\,.\eqno(14)$$ (14) and the definition of $Q_*$ imply: $$Q(\chi
_{\omega '}\chi _{\omega ''}) = \chi _{\omega '}(\varphi ) Q(\chi
_{\omega ''}),~~\omega '\in \Sigma _1,~ \omega ''\in (\Sigma
_1)_0.\eqno(15)$$ Equalities (15) permit to define the unique
measurable function $q$ on $\Omega _2$ such that $$Q(\chi _\sigma ) =
q\chi _\sigma (\varphi ),~~ \sigma \in (\Sigma _1)_0\eqno(16)$$ By
Theorem 2 of [Za3], (16) implies (12). Since the operator $Q$ is
inversible, the transformation $\varphi $ is also inversible. The
theorem is proved.

\rque Comments|

1. The theorem is true for more general measure spaces. In
particular, the proof did not use continuity of the measure.

2. To prove Theorem 1 in the real case, probably it is impossible
simply to pass to the complexifications. Indeed, in general there is
no universal definition of a norm in the complexification of an ideal
function space that would be at the same time ideal and hold the
property of ``extension of isometries'', \ie such that all of them
extend to isometries of the complexification. The simplest example of
such situation is the Minkowski plane $E^2$, where the unit sphere is
a regular octagon. This symmetric space has an extra isometry, namely
the rotation on angle $\pi /4$~\footnote{(***)}{This example was also
noted by Yu. Sokolovski.}. It is proved in [BrSe] that extra
isometries of real symmetric sequence spaces could exist in dimensions
2 and 4 only. Due to Theorem 3 from [Ta] any complex symmetric
sequence space permits just standard isometries, i.e. permutations and
rotations of coordinates.

3. It would be interesting to extend Theorem 1 to more general classes of
Banach lattices; for instance, to describe the class of Banach
lattices in which any invertible isometry is disjoint, \ie maps
disjoint elements to disjoint elements.


\titre Literature|

{\ept

\livre AbZa|Yu. Abramovich, M. Zaidenberg|A rearrangement invariant
space isometric to $L_p$ coincides with $L_p$|in: N. Kalton and
E. Saab eds.,
 Proceedings of the Conference {\ \it Interaction between Functional
Analysis, Harmonic Analysis, and Probability}, Marcel Dekker|{\rm ~New
York (forthcoming)}|

\article BrSe|M.Sh. Braverfman, E.M. Semenov|Isometries of symmetric
spaces|Soviet Math. Dokl.|15|1974|1027--1030|

\article FlJe|R.J. Fleming, J.E. Jamison|Isometries on Banach spaces
-- a survey|preprint||1994|1--48|

\article KaRa1|N.J. Kalton, B. Randrianantoanina|Isometries on
rearrangement invariant
spaces|C. R. Acad. Sci. Paris|316|1993|351--355|

\article KaRa2|N.J. Kalton, B. Randrianantoanina|Surjective isometries
on rearrangement invariant spaces|Quart.J.Math.Oxford|(2),
45|1994|301--327|

\article KaWo|N.J. Kalton, G.W. Wood|Orthonormal systems in Banach
spaces and their
applications|Math. Proc. Cambr. Phil. Soc.|79|1976|493--510|

\article PKL|Pei-Kee Lin|Elementary isometries of rearrangement
invariant spaces|preprint||1995|1--24|

\article Lu|G. Lumer|On the isometries of reflexive Orlicz
spaces|Ann. Inst. Fourier, Grenoble|13|1963|99--109|

\article Ra|B. Randrianantoanina|Isometric classification of norms in
rearrangement - invariant function spaces|preprint||1995|1--15|

\livre Ro|S. Rolewicz|Metric Linear Spaces|Polish Scientific
Publishers, Warsaw|1972|

\livre Se|E.M. Semenov|Doctorate dissertation|University of
Voronezh|1963 {\rm (in Russian)}|

\livre SMB|Functional analysis|Edited by S.G. Krein, Nauka|Moscow|1972
{\rm (in Russian)}|

\article SkZa1|A.I. Skorik, M.G. Zaidenberg|Groups of isometries
containing reflexions|Function. Anal. Appl.|10|1976|322--323|

\article SkZa2|A. Skorik, M. Zaidenberg|On isometric reflections in
Banach spaces|Preprint, Institut Fourier, Laboratoire de
Math\'ematiques, Grenoble|267|1994|1--36|

\article Ta|K.W. Tam|Isometries of certain function spaces|Pacific
J. Math.|31|1969|233--246|

\article To|G. Tom\u si\u c|Skew-symmetric operators on real Banach
spaces|Glasnik Mat.|12|1977|93--97|

\article Zab|P.P. Zabreiko|Ideal function spaces|Vestnik of the
University of Yaroslavl'|8|1974|12--52 (in Russian)|

\article Za1|Zaidenberg M.G|On the isometric classification of
symmetric spaces|Soviet Math. Dokl.|18|1977|636--640|

\article Za2|M.G. Zaidenberg|Special representations of isometries of
functional spaces|in: Investigations on the theory of functions of
several real variables|{\rm Edited by Yu.A. Brudnyi,
Yaroslavl'}|1980|84--91 (in Russian)|

\article Za3|M.G. Zaidenberg|Groups of isometries of Orlich
spaces|Soviet Math. Dokl.|17|1976|432--436|

}

\iffin
\end